\documentclass[11pt]{article}
\usepackage{amssymb,latexsym,amscd,amsmath,amsfonts,enumerate,supertabular}
\usepackage{xfrac}
\usepackage{breqn}
\usepackage{graphicx}
\usepackage{tabularx}
\usepackage{caption}
\usepackage{epstopdf}
\usepackage[top=2cm, bottom=2cm, left=2.5cm, right=2cm]{geometry}
\usepackage{amssymb,latexsym,amscd,amsmath,amsfonts,enumerate,supertabular}
\usepackage{graphicx}
\usepackage{epstopdf}
\usepackage{tabularx}
\usepackage{caption}

\newtheorem{Theorem}{Theorem}[section]

\newtheorem{Proposition}[Theorem]{Proposition}

\begin{document}

\title{\bf Existence results and numerical solution of fully fourth order nonlinear functional differential equations }

\author{ Dang Quang A$^{\text a}$,  Nguyen Thanh Huong$^{\text b}$, Dang Quang Long$^{\text c}$\\
\small $^{\text a}$ { Center for Informatics and Computing, Vietnam Academy of Science and Technology}\\
\small{  18 Hoang Quoc Viet, Cau Giay, Ha Noi, Viet Nam}\\
\small{ Email address: dangquanga@cic.vast.vn}\\
\small $^{\text b}$ { College of Sciences, Thai Nguyen University, Thai Nguyen, Viet Nam}\\
\small { Email address: nguyenthanhhuong2806@gmail.com}\\
$^{\text c}$ {\small Institute of Information Technology, VAST,}\\
{\small 18 Hoang Quoc Viet, Cau Giay, Hanoi, Vietnam}\\
{\small Email: dqlong88@gmail.com}}
\date{ }
\maketitle
\noindent {\bf Abstract. }
In this paper we consider a  boundary value problem for fully fourth order nonlinear functional differential equation which contains all lower derivatives of proportional delay arguments.
By the reduction of the problem to operator equation for the right hand side nonlinear function we establish the existence and uniqueness of solution and construct iterative  methods on both continuous and discrete levels for solving it. We obtain the total error estimate for the discrete iterative solution.  Many examples demonstrate the validity of the obtained theoretical results and the efficiency of the numerical method. \par

\section{Introduction}
In this paper we consider a class of boundary value problems for fourth order nonlinear functional differential equation
\begin{equation}\label{eq1}
\begin{split}
&u''''(t)=f(t,U(t)), \quad 0<t<1,\\
&u(0)=a, u(1)= b, u'(0)=c, u'(1)=d,
\end{split}
\end{equation}
where
\begin{equation}\label{eq2}
U(t)=(u(t), u(\varphi_0(t)), u'(t), u'(\varphi_1(t)), u''(t), u''(\varphi_2(t)), u'''(t), 
u'''(\varphi_3(t)) ),
\end{equation}
and  $f: [0, 1] \times \mathbb{R}^8 \rightarrow \mathbb{R}$ and  $\varphi_i: [0, 1] \rightarrow [0, 1]\; (i=\overline{0,3})$ are continuous functions.\par
A particular case of the problem \eqref{eq1} when $U(t)=(u(t), u(\varphi(t))$ was studied in \cite{Bica}. In this work by reducing the problem to the equivalent integral equation the authors established the existence and uniqueness of a solution under very strong assumptions, one of them is the Lipschitz condition of the function $f(t,u,v)$ in the variables $u, v$ in the domain $[a,b] \times \mathbb{R}^2$. Under some more other assumptions the authors constructed iterative scheme with the use of cubic spline interpolation at each iteration. The analysis of convergence was made, but it is a regret that in the proof of Theorem 8 \cite[page 139]{Bica} there  are vital errors shown by us in  a private letter to the authors and in \cite{Dang1}. After that these errors were overcome in the  corrigendum \cite{Bica1}. It should be emphasized that in all illustrative examples (Example 11 and Example 12 \cite{Bica}) the conditions for ensuring the existence and uniqueness of solutions are not satisfied although the numerical results show good convergence.\par
In a very recent work \cite{Dang1} we proposed an effective method for investigating simultaneously the existence of solution and numerical method for the third order functional differential equation
\begin{equation}\label{eqAPNUM}
u''' (t)=f(t,u(t),u(\varphi(t)))
\end{equation}
associated with general linear boundary conditions. The method is a further development of our method for nonlinear boundary value problems in, e.g., \cite{Dang3}, \cite{Dang4}.\par 
In connection to the general functional differential equation $u^{(p)}(t)= f(t,u(\varphi (t)), p \ge 1$ it should be said that an its particular case when $\varphi (t) = \alpha t, 0< \alpha <1$, the equation is called pantorgraph equation or proportional delay differential equation. This type equation is intensively investigate in many works such as \cite{Chamekh}, \cite{Jafari},  \cite{Raja}, \cite{Reut}, \cite{Suay}. The above mentioned works are concerned of initial value problems.\par 
 The solution of boundary value problems for second order functional differential equation is considered in several works, e.g., \cite{Bica0}, \cite{Hou}, \cite{Khuri}, \cite{Wazwaz},  where numerical, semi-analytical  and neural network methods are used.\par
The aim of the present work is to investigate the existence and uniqueness of solution of the boundary value problem for general  fourth order nonlinear functional differential equation \eqref{eq1}-\eqref{eq2} and develop an efficient numerical method for finding the solution.

\section{Existence and uniqueness of solution}
Following the approach in \cite{Dang1,Dang2} (see also \cite{Dang3,Dang4}) to investigate the problem \eqref{eq1}-\eqref{eq2} we introduce the nonlinear operator $A$ defined in the space of continuous functions $C[0,1]$ by the formula:
\begin{equation}\label{eq3}
(A\psi )(t) =f(t,U(t)),
\end{equation}
where $u(t)$ is the solution of the problem
\begin{equation}\label{eq4}
\begin{split}
u'''(t)&=\psi (t), \quad 0<t<1\\
u(0)=a, u(1)&= b, u'(0)=c, u'(1)=d.
\end{split}
\end{equation}
It is easy to verify the following
\begin{Proposition}\label{Prop1}
If the function $\psi$ is a fixed point of the operator $A$, i.e., $\psi$ is the solution of the operator equation 
\begin{equation}\label{eq5}
A\psi = \psi ,
\end{equation} 
where $A$ is defined by \eqref{eq3}-\eqref{eq4}
then the function $u(t)$ determined from the BVP \eqref{eq4} is a solution of the BVP \eqref{eq1}-\eqref{eq2}. Conversely, if the function $u(x)$ is the solution of the BVP \eqref{eq1}-\eqref{eq2} then the function 
\begin{equation*}
\psi (t)= f(t, U(t) )
\end{equation*}
satisfies the operator equation \eqref{eq5}.
\end{Proposition}
Now, let $G(t,s) $ be the Green function of the problem \eqref{eq4}. We have
\begin{equation*}
\begin{aligned}
G(t,s)=\dfrac{1}{6}\left\{\begin{array}{ll}
s^2(1-t)^2(3t-s-2ts), \quad 0\le s \le t \le 1,\\
t^2(1-s)^2,(3s-t-2ts) \quad 0\le t \le s \le 1,\\
\end{array}\right.
\end{aligned}
\end{equation*}
Denote $G_0(t,s)=G(t,s)$ and
\begin{equation}\label{eq6}
G_i(t,s)=\frac{\partial^iG(t,s)}{\partial t^i},\quad (i=\overline{1,3}).
\end{equation}
Then there hold the estimates
\begin{equation}\label{eq7}
\int_0^1 |G_i(t,s)| ds \le M_i \quad (i=\overline{0,3})
\end{equation}
with
\begin{equation}\label{eq8}
M_0=\frac{1}{384}, M_1=\frac{1}{72\sqrt{3}}, M_2=\frac{1}{12}, M_3=\frac{1}{2}.
\end{equation}
The solution of the problem \eqref{eq4} is represented in the form
\begin{equation}\label{eq9}
u(t)= g(t)+ \int_0^1 G(t,s)\psi (s) ds,
\end{equation}
where $g(t) $ is the polynomial of second degree satisfying the boundary conditions
\begin{equation}\label{eq10}
g(0)=a, g(1)= b, g'(0)=c, g'(1)=d.
\end{equation}
Taking the derivatives of \eqref{eq8} we obtain
\begin{equation}\label{eq11}
u^{(i)}(t)= g^{(i)}(t)+ \int_0^1 G_i(t,s)\psi (s) ds, \; (i=\overline{0,3})
\end{equation}
Therefore, from \eqref{eq6} we have
\begin{equation}\label{eq12}
\|u^{(i)}\| \le \|g^{(i)}\| +M_i \|\psi\|,\; (i=\overline{0,3}),
\end{equation}
where $\|.\| $ is the norm in the space $C[0,1]$, $\|v\|= \max _{0\le t\le 1}|v(t)|, v\in C[0,1]$.\\
Now for every positive number $M$ define the domain
\begin{equation}\label{eq13}
\begin{split}
\mathcal{D}_M=\Big \{ (t,u,\bar{u},y,\bar{y},v,\bar{v},z,\bar{z}) \mid  0\leq t\leq 1;\; &|u|,|\bar{u}|\leq \|g\|+M_0 M; |y|,|\bar{y}| \leq \|g'\|+M_1 M,  \\
 &|v|, |\bar{v}| \le \|g''\|+M_2 M;  |z|, |\bar{z}| \le \|g'''\|+M_3 M \Big \},
 \end{split}
\end{equation}
As usual, we denote by $B[0,M]$ the closed ball of the radius $M$ centered at $0$ in the space of continuous functions $C[0,1]$. 
\begin{Theorem}\label{thm1}
Assume that:
\begin{description}
\item (i) The function $\varphi_k (t)\; (k=\overline{0,3})$ are continuous functions from $[0,1]$ to $[0,1]$.
\item (ii) The function $f(t,u,\bar{u},y,\bar{y},v,\bar{v},z,\bar{z})$ is continuous and bounded by $M$ in the domain $\mathcal{D}_M$, i.e.,
\begin{equation}\label{eq14}
|f(t,u,\bar{u},y,\bar{y},v,\bar{v},z,\bar{z})|\le M \quad \forall (t,u,\bar{u},y,\bar{y},v,\bar{v},z,\bar{z}) \in \mathcal{D}_M.
\end{equation}
\item (iii) The function $f(t,u,\bar{u},y,\bar{y},v,\bar{v},z,\bar{z})$ satisfies the Lipschitz conditions in the variables $u,\bar{u},y,\bar{y}$, $v,\bar{v},z,\bar{z}$ with the coefficients $L_i \ge 0, \; (i= \overline{0,7})$ in $\mathcal{D}_M$, i.e.,
\begin{equation}\label{eq15}
\begin{split}
|f(t,u_2,\bar{u}_2,y_2,\bar{y}_2,v_2,\bar{v}_2,z_2,\bar{z}_2)-f(t,u_1,\bar{u}_1,y_1,\bar{y}_1,v_1,\bar{v}_1,z_1,\bar{z}_1) |\le L_0 |u_2-u_1|+L_1 |\bar{u}_2-\bar{u}_1| \quad \\
+L_2 |y_2-y_1|+L_3 |\bar{y}_2-\bar{y}_1|+ L_4 |v_2-v_1|+L_5 |\bar{v}_2-\bar{v}_1|+L_6 |z_2-z_1|+L_7 |\bar{z}_2-\bar{z}_1| \\
\forall (t,u_i,\bar{u}_i,y_i,\bar{y}_i,v_i,\bar{v}_i,z_i,\bar{z}_i)\in \mathcal{D}_M\;
(i=1,2)
\end{split}
\end{equation}
\item (iv)   
\begin{equation}\label{eq16}
q:= (L_0+L_1)M_0 +(L_2+L_3)M_1+(L_4+L_5)M_2+(L_6+L_7)M_3 <1.
\end{equation}

\end{description}
The the problem \eqref{eq1}-\eqref{eq2} has a unique solution $u(t) \in C^4[0,1]$, satisfying the estimate
\begin{equation}\label{eq17}
|u^{(i)}(t)| \le \|g^{(i)}\| +M_iM , (i=\overline{0,3})   \quad \forall t\in[0,1].
\end{equation}
\end{Theorem}
\noindent {\bf Proof.} Sketch the steps of the proof:\\
\noindent i) Show that the operator $A$ is a mapping  $B[0,M]\rightarrow B[0,M]$.\\
\noindent ii)Show that $A$ is a contraction mapping in $B[0,M]$.\\
Then the operator equation \eqref{eq5} has a unique solution $\psi \in B[0,M]$. By Proposition \ref{Prop1} the solution of the problem \eqref{eq4} for this right-hand side $\psi (t)$ is the solution of the original problem \eqref{eq1}-\eqref{eq2}.\\
The detailed proof is similar to that in \cite{Dang1} for the equation \eqref{eqAPNUM}.

\section{Solution method and its convergence}\label{sec3}
Consider the following iterative method:
\begin{enumerate}
\item Given $\psi_0 \in B[0,M]$, for example,
\begin{equation}\label{iter1c}
\psi_0(t)=f(t,0,0,0,0,0,0,0,0).
\end{equation}
\item Knowing $\psi_k(t)$  $(k=0,1,...)$ compute
\begin{equation}\label{iter2c}
\begin{split}
u_k(t) &= g(t)+\int_0^1 G(t,s)\psi_k(s)ds ,\\
\bar{u}_k(t) &= g(\varphi_0 (t))+\int_0^1 G(\varphi_0(t),s)\psi_k(s)ds \\
y_k(t) &= g'(t)+\int_0^1 G_1(t,s)\psi_k(s)ds ,\\
\bar{y}_k(t) &= g'(\varphi_1 (t))+\int_0^1 G_1(\varphi_1(t),s)\psi_k(s)ds \\
v_k(t) &= g''(t)+\int_0^1 G_2(t,s)\psi_k(s)ds ,\\
\bar{v}_k(t) &= g''(\varphi_2 (t))+\int_0^1 G_2(\varphi_2(t),s)\psi_k(s)ds \\
z_k(t) &= g'''(t)+\int_0^1 G_3(t,s)\psi_k(s)ds ,\\
\bar{z}_k(t) &= g'''(\varphi_3 (t))+\int_0^1 G_3(\varphi_3(t),s)\psi_k(s)ds .\\
\end{split}
\end{equation}
\item Update
\begin{equation}\label{iter3c}
\psi_{k+1}(t) = f(t,u_k(t),\bar{u}_k(t),y_k(t),\bar{y}_k(t),v_k(t),\bar{v}_k(t),z_k(t),\bar{z}_k(t) ).
\end{equation}
\end{enumerate}

Set 

\begin{equation}\label{eqpkd}
p_k=\dfrac{q^k}{1-q} ,\; d=\|\psi _1 -\psi _0\|.
\end{equation} 
\begin{Theorem}[Convergence]\label{thm2} Under the assumptions of Theorem \ref{thm1} the above iterative method converges and there holds the estimate
\begin{equation*}
\|u^{(i)}_k-u^{(i)}\| \leq M_ip_kd, \; (i=\overline{0,3}),
\end{equation*} 
where $u$ is the exact solution of the problem \eqref{eq1}-\eqref{eq2}, $u^{(i)}$ is its derivative of order $i$ and $M_i$ is given by \eqref{eq8}.
\end{Theorem}
This theorem follows straightforward from the convergence of the iterative method for  fixed point of the operator $A$, the representations \eqref{eq9}, \eqref{eq11} and the formulas for computing $u_k(t),y_k(t),v_k(t),z_k(t) $ in \eqref{iter2c}.\par

Now we design a discrete iterative scheme for realizing the above iterative method on continuous level.
 For this purpose, we construct the uniform grid $\bar{\omega}_h=\{t_i=ih, \; h=a/N, i=0,1,...,N  \}$ on the interval $[0, 1]$  
  and denote by $\Phi_k(t), U_k(t),\bar{U}_k(t), Y_k(t),\bar{Y}_k(t)$,  $V_k(t),\bar{V}_k(t)$,$Z_k(t),\bar{Z}_k(t)$ the grid functions, which are defined on the grid $\bar{\omega}_h$ and approximate the functions $\psi_k (t), u_k(t), \bar{u}_k(t),  y_k(t), \bar{y}_k(t)$,  $v_k(t), \bar{v}_k(t), z_k(t), \bar{z}_k(t) $ on this grid, respectively.\par
  Below are steps of the discrete iterative method:\\
\begin{enumerate}
\item Given 
\begin{equation}\label{iter1d}
\Psi_0(t_i)=f(t_i,0,0,0,0,0,0),\ i=0,...,N. 
\end{equation}
\item Knowing $\Psi_k(t_i),\; k=0,1,...; \; i=0,...,N, $  compute approximately the definite integrals \eqref{iter2c} by the trapezoidal rule
\begin{equation}\label{iter2d}
\begin{split}
U_k(t_i) &=g(t_i) +\sum _{j=0}^N h\rho_j G(t_i,t_j)\Psi_k(t_j), \\
\bar{U}_k(t_i) &= g(\xi_{0i})+\sum _{j=0}^N h\rho_j G(\xi_{0i},t_j)\Psi_k(t_j) , \\
Y_k(t_i) &=g'(t_i) +\sum _{j=0}^N h\rho_j G_1(t_i,t_j)\Psi_k(t_j), \\
\bar{Y}_k(t_i) &= g'(\xi_{1i})+\sum _{j=0}^N h\rho_j G_1(\xi_{1i},t_j)\Psi_k(t_j) , \\
V_k(t_i) &=g''(t_i) +\sum _{j=0}^N h\rho_j G_2(t_i,t_j)\Psi_k(t_j), \\
\bar{V}_k(t_i) &= g''(\xi_{2i})+\sum _{j=0}^N h\rho_j G_2(\xi_{2i},t_j)\Psi_k(t_j) , \\
Z_k(t_i) &=g'''(t_i) +\sum _{j=0}^N h\rho_j G_3^{*}(t_i,t_j)\Psi_k(t_j), \\
\bar{Z}_k(t_i) &= g'''(\xi_{3i})+\sum _{j=0}^N h\rho_j G_3^{*}(\xi_{3i},t_j)\Psi_k(t_j) , \  (i=0,...,N),\\
\end{split}
\end{equation}
\noindent where $\rho_j$ are the weights 
\begin{equation*}
\rho_j = 
\begin{cases}
1/2,\; j=0,N\\
1, \; j=1,2,...,N-1
\end{cases}, \quad \xi_{mi}=\varphi_m (t_i),\; m=\overline{0,3}
\end{equation*}
and 
\begin{equation*}
G_3^{*}(t,s) = 
\begin{cases}
G_{3}(t,s) , \; s \neq t, \\
\frac{1}{2}[ \lim _{s\rightarrow t-0}G_{3}(t,s)+ \lim _{s\rightarrow t+0}G_{3}(t,s)], \; s=t.
\end{cases} 
\end{equation*}

\item Update
\begin{equation}\label{iter3d}
\Psi_{k+1}(t_i) = f(t_i,U_k(t_i),\bar{U}_k(t_i), Y_k(t_i),\bar{Y}_k(t_i),V_k(t_i),\bar{V}_k(t_i),Z_k(t_i),\bar{Z}_k(t_i)).
\end{equation}
\end{enumerate}
To study the convergence of the above discrete iterative method we need some auxiliary results.
\begin{Proposition}\label{prop2}
Assume that the functions $\varphi_m(t), \; (m=\overline{0, 3})$ have continuous derivatives up to second order in $[0,1]$ and the function $f(t,u,\bar{u},y,\bar{y},v,\bar{v},z,\bar{z})$ has all 
  partial derivatives continuous up to second order in the domain $\mathcal{D}_M$. 
   Then the functions generated  by the iterative method \eqref{iter1c}-\eqref{iter3c} satisfy the following conditions:
   \begin{equation}\label{eqSmooth}
   \begin{aligned}
 &z_k(t) \in C^3[0,1], v_k(t) \in C^4[0,1], y_k(t) \in C^5[0,1], u_k(t) \in C^6[0,1],\\
 & \psi_k(t),  \bar{u}_k(t),  \bar{y}_k(t),  \bar{v}_k(t),  \bar{z}_k(t)\in C^2[0,1].
\end{aligned}
   \end{equation}

\end{Proposition}
\noindent {\bf Proof.} We shall prove the above proposition by induction for $k$.\\
For $k=0$, by the assumption on the function $f$ we have $\psi_0(t)=f(t,0,0,0,0,0,0) \in C^2[0,1]$. Taking into account the expression of $G_3(t,s)$ (see Appendix ) we have
\begin{align*}
z_0(t)&= g'''(t) +\int_0^1 G_3(t,s) \psi_0(s) ds \\
&= g'''(t) +\int_0^1 s^2(3-2s)\psi_0(s) ds  - \int_t^1 \psi_0(s) ds.
\end{align*}
Therefore, $z_0'(t)=\psi_0(t)$. It implies that $z_0'(t)\in C^2[0,1] $ and consequently, $z_0(t)\in C^3[0,1] $. In view of \eqref{iter2c} it follows that $v_0(t)\in C^4[0,1], y_0(t)\in C^5[0,1], u_0(t)\in C^6[0,1] $. Next, we can obtain $\bar{z}_0'(t)= \varphi_3'(t) \psi_0 (\varphi_3(t))\in C^1[0, 1]$. Hence, $\bar{z}_0(t)\in C^2[0,1]$. It is also easy to show that $\bar{v}_0(t), \bar{y}_0(t), \bar{u}_0(t)\in C^2[0,1] $. \\
Thus, the proposition is true for $k=0$.\par 
Now suppose that the proposition is true for $k \ge 0$, i.e., there holds \eqref{eqSmooth}. Then from the assumption of the smoothness of the function $f$ and from the formula  \eqref{iter3c} it follows that $\psi_{k+1}(t) \in C^2[0, 1]$. Making similar argument as for the case $k=0$ we obtain
  \begin{equation*}
   \begin{aligned}
 &z_{k+1}(t) \in C^3[0,1], v_{k+1}(t) \in C^4[0,1], y_{k+1}(t) \in C^5[0,1], u_{k+1}(t) \in C^6[0,1],\\
 &  \bar{u}_{k+1}(t),  \bar{y}_{k+1}(t),  \bar{v}_{k+1}(t),  \bar{z}_{k+1}(t)\in C^2[0,1].
\end{aligned}
   \end{equation*}
Thus, the proposition is proved.
\begin{Proposition}\label{prop3}
For any function $\psi (t) \in C^2[0,1]$ there hold the estimates
\begin{align}
\int_0^1 G_m (t_i,s) \psi (s) ds = \sum _{j=0}^N h\rho_j G_m(t_i,s_j)\psi(s_j) +O(h^2), \; m=0,1,2, \label{eqTp1}\\ 
\int_0^1 G_3 (t_i,s) \psi (s) ds = \sum _{j=0}^N h\rho_j G_3^{*}(t_i,s_j)\psi(s_j) +O(h^2),  \label{eqTp2}\\ 
\int_0^1 G_m (\xi_{mi},s) \psi (s) ds = \sum _{j=0}^N h\rho_j G_m(\xi_{mi},s_j)\psi(s_j) +O(h^2),\; m=0,1,2, \label{eqTp3}\\
\int_0^1 G_3 (\xi_{3i},s) \psi (s) ds = \sum _{j=0}^N h\rho_j G_3^{*}(\xi_{3i},s_j)\psi(s_j) +O(h^2),\label{eqTp4}
\end{align}
where $\xi_{mi}=\varphi_m (t_i),\; m=\overline{0,3}$.
\end{Proposition}
\noindent {\bf Proof} For proof of \eqref{eqTp1}- \eqref{eqTp2} see \cite[Proposition 2]{Dang5}, and for \eqref{eqTp3}-\eqref{eqTp4} see \cite[Proposition 3.3]{Dang1}.

\begin{Proposition}\label{prop4}
Under the assumptions of Theorem \ref{thm1} we have the estimates
\begin{align}
\| \Psi _k -\psi_k\|_{\bar{\omega_h}} =O(h^2), \; \| U _k -u_k\|_{\bar{\omega_h}} =O(h^2), \label{est1}\\
\| Y _k -y_k\|_{\bar{\omega_h}} =O(h^2), \| V _k -v_k\|_{\bar{\omega_h}} =O(h^2),
\| Z _k -z_k\|_{\bar{\omega_h}} =O(h^2),\label{est2}
\end{align}
where $\|. \|_{\bar{\omega_h}}$ is the max-norm of grid function defined on the grid $\bar{\omega_h}$.
\end{Proposition}
\noindent {\it Proof.} It is not hard to prove the proposition by induction using Propositions \ref{prop2} and \ref{prop3}, and taking into account the formulas \eqref{iter1c}, \eqref{iter1d} and \eqref{iter2c}, \eqref{iter2d} of the iterative methods on continuous and on discrete levels. \par 
Similar but simpler propositions are proved in detail in \cite[Proposition 3.5]{Dang1} and \cite[Proposition 3]{Dang5}.\par

Now combining  Proposition \ref{prop4} with Theorem \ref{thm2} we obtain the following result on the error estimate of the actually obtained numerical solution of the original problem.
\begin{Theorem}\label{thm3}
Under the assumptions of Theorem \ref{thm1} for the approximate solution of the problem \eqref{eq1}-\eqref{eq2} obtained by the discrete iterative method \eqref{iter1d}-\eqref{iter3d} we have the estimates
\begin{equation*}
\begin{split}
\|U_k-u\|_{\bar{\omega_h}} \leq M_0p_kd +O(h^2), \; \|Y_k-u'\|_{\bar{\omega_h}} \leq M_1p_kd +O(h^2), \\
\|V_k-u''\|_{\bar{\omega_h}} \leq M_2p_kd +O(h^2), \; \|Z_k-u'''\|_{\bar{\omega_h}} \leq M_3p_kd +O(h^2), \\
\end{split}
\end{equation*}
where $p_k$ and $d$ are defined by \eqref{eqpkd}.
\end{Theorem}

\section{Examples}
In all  examples below first we verify the conditions of Theorem \ref{thm1} which guaranty the existence and uniqueness of a solution and also the convergence of the iterative methods. After that 
 we perform the iterative method \eqref{iter1d}-\eqref{iter3d} until $\| U _k -U_{k-1}\|_{\bar{\omega_h}} \le 10^{-16}$. 
 In the tables of results for the convergence of the iterative method $N+1$ is the number of grid points, $h=1/N$, 
 $K$ is the number of iterations performed, $u$ is the exact solution if it is known and  
 $Error=\|U_K-u \|_{\bar{\omega_h}}$ and $Error1=\|Y_K-u' \|_{\bar{\omega_h}}$ are errors of the solution and its first derivative, respectively.
 
\noindent {\bf Example 1 } \cite[Example 11]{Bica}. Consider the following problem
\begin{equation}\label{eqExam1}
\begin{split}
u^{(4)}(t)= \frac{22}{(t+1)^5}+\frac{1}{(t+1)^2}\Big ( [u(t)]^2+[u(t)]^3   \Big )u(\frac{t}{2}),\; 0<t<1,\\
u(0)=1, \ u(1)= \frac{1}{2}, \ u'(0)=-1, \ u'(1)=-\frac{1}{4}
\end{split}
\end{equation}
for which the exact solution is $u(t)=\frac{1}{t+1}$.\\
For this problem the right hand side function $f$ is
\begin{equation*}
f(t,u,\bar{u})=\frac{22}{(t+1)^5}+\frac{1}{(t+1)^2} (u^2 + u^3)\bar{u}.
\end{equation*}
Obviously, this function does not satisfy the Lipschitz conditions in the domain $[0,1] \times \mathbb{R}^2$.
Thus, the important condition (ii) in \cite[page 131]{Bica} is not met. So, the uniqueness of the solution is not ensured, and of course, Theorem 2 and Theorem 8 there on the convergence of the iterative method are not applicable. Differently from Bica \cite{Bica} below we show that the problem \eqref{eqExam1} has a unique solution and the discrete iterative method \eqref{iter1d}-\eqref{iter3d} converges with accuracy of second order. Indeed, for the problem it is easy to find the function $g(t)=-\frac{1}{4}t^3+\frac{3}{4}t^2-t+1$ satisfying the boundary conditions, and $\|g\|=\|g'\| =1, \|g''\|=\|g'''\| = \frac{3}{2}$. It is possible to verify that for $M=25$ we have
$|f(t,u,\bar{u}| \le M $ in the domain $\mathcal{D}_M=\Big \{ (t,u,\bar{u}) \mid  0\leq t\leq 1;\; |u|,|\bar{u}|\leq 1+M_0 M \}$. Further, in $\mathcal{D}_M$ the function $f(t,u,\bar{u})$ satisfies the Lipschitz conditions in $u$ and $\bar{u}$ with the coefficients $L_0=6, L_1=2.4$. Therefore, $q= (L_0+L_1)M_0 =0.0219<1$. By Theorem \ref{thm1} the problem \eqref{eqExam1} has a unique solution and by Theorem \ref{thm3} the iterative method \eqref{iter1d}-\eqref{iter3d} converges. \par
Below, in Table \ref{table1} we report  the convergence of the iterative method . 
\begin{table}[ht!]
\centering
\caption{The convergence in Example 1. }
\label{table1}
\begin{tabular}{ccccc}
\hline 
$N$ &	$h^2$ &	$K$	&$Error$&$Error1$  \\
\hline 
50	& 4.0000e-04 &	8&	1.4520e-08 & 3.8031e-08\\
100	& 1.0000e-04 &	8&	9.0870e-10 & 2.3801e-09\\
150	&4.4444e-05	&  8&1.7954e-10&  4.7028e-10 \\
200	&2.5000e-05	&8	&5.6812e-11& 1.4881e-10\\
300	&1.1111e-05	&9	&1.1223e-11 & 2.9396e-11\\
400&	6.2500e-06&	9	&3.5512e-12& 9.3012e-12\\
500	&4.0000e-06	&8	&1.4546e-12& 3.8097e-12\\
800	&1.5625e-06&	9&	 2.2204e-13&  5.8115e-13\\
1000&	1.0000e-06	&9	&9.1038e-14& 2.3798e-13\\
\hline 
\end{tabular} 
\end{table}
From the table we see that the errors of the approximate solution and its first derivative are of order 4 although Theorem \ref{thm3} ensures only  order 2. This fact also will be observed in the next example.
\\

\noindent {\bf Example 2 } \cite[Example 12]{Bica}. Consider the following problem
\begin{equation}\label{eqExam2}
\begin{split}
u^{(4)}(t)= e^{-t} [u(t)]^{3/2} u(\frac{1}{2}),\; t\in (0,1)\\
u(0)=1, \ u(1)= e, \ u'(0)=1, \ u'(1)=e
\end{split}
\end{equation}
for which the exact solution is $u(t)=e^{t}$.\\
For this problem the right hand side function $f$ is
\begin{equation*}
f(t,u,\bar{u})=e^{-t} u^{3/2}\bar{u}.
\end{equation*}
As in the previous example, the function $f(t,u,\bar{u})$ does not satisfy the Lipschitz conditions in the domain $[0,1] \times \mathbb{R}^2$. So, 
the important condition (ii) in \cite[page 131]{Bica} is not met. Hence, the uniqueness of the solution is not ensured, and of course, Theorem 2 and Theorem 8 there on the convergence of the iterative method are not applicable.
Below we apply the theory in the previous section to the example. We have for this example
$ g(t)=(3-e)t^3 +(2e-5)t^2 +t+1$, and consequently, $\|g\|=e$.
It is possible to verify that for $M=15$ we have $|f(t,u,\bar{u}| \le M $ in the domain $\mathcal{D}_M=\Big \{ (t,u,\bar{u}) \mid  0\leq t\leq 1;\; |u|,|\bar{u}|\leq e+M_0 M \}$. Further, in $\mathcal{D}_M$ the function $f(t,u,\bar{u})$ satisfies the Lipschitz conditions in $u$ and $\bar{u}$ with the coefficients $L_0=7, L_1=5$. Therefore, $q= (L_0+L_1)M_0 =0.0313 <1$. By Theorem \ref{thm1} the problem \eqref{eqExam1} has a unique solution and by Theorem \ref{thm3} the iterative method \eqref{iter1d}-\eqref{iter3d} converges. \par
The results of convergence  of the iterative method for this example is given in Table \ref{table2}.
\begin{table}[ht!]
\centering
\caption{The convergence in Example 2. }
\label{table2}
\begin{tabular}{ccccc}
\hline 
$N$ &	$h^2$ &	$K$	&$Error$&$Error1$  \\
\hline 
50	& 4.0000e-04 &	9&	3.0102e-10 & 2.2923e-09\\
100	& 1.0000e-04 &	9&	1.8814e-11 & 1.4328e-10\\
150	&4.4444e-05	&  9&3.7157e-12 &  2.8302e-11 \\
200	&2.5000e-05	&10	&1.1755e-12& 8.9548e-12 \\
300	&1.1111e-05	& 10	&2.3226e-13  & 1.7686e-12 \\
400&	6.2500e-06&	10	&7.3275e-14& 5.5933e-13\\
500	&4.0000e-06	& 9	& 3.0198e-14& 2.2893e-13\\
800	&1.5625e-06&	9&	 4.8850e-15 &  3.4639e-14\\
1000&	1.0000e-06	&9	&2.2204e-15& 1.4211e-14\\
\hline 
\end{tabular} 
\end{table}

\noindent {\bf Example 3 .} Consider a more complicated problem
\begin{equation}\label{eqExam3}
\begin{split}
u^{(4)}(t)&= e^{t}+\frac{1}{9}\Big( [u(\frac{t}{2})]^2 u'''(\frac{t}{2})-u'(t)u''(\frac{t}{2})+u(t)u'''(t)-[u''(t)]^2   \Big),\; t\in (0,1)\\
u(0)&=1, \ u(1)= e, \ u'(0)=1, \ u'(1)=e
\end{split}
\end{equation}
for which the exact solution is $u(t)=e^{t}$.\\
For this example 
\begin{align*}
&f=f(t,u,\bar{u},y,\bar{y},v,\bar{v},z,\bar{z})= e^{t}+\frac{1}{9}(\bar{u}^2 \bar{z}-y\bar{v}+uz-v^2),\\
&g(t)=(3-e)t^3 +(2e-5)t^2 +t+1.
\end{align*}
So, we have $\|g\|=\|g'\|=e, \|g''\|=8-2e, \|g'''\| =18-6e$. It is possible to check that for $M=20$ there holds
$|f(t,u,\bar{u},y,\bar{y},v,\bar{v},z,\bar{z})| \le M$ in $\mathcal{D}_M$, where $\mathcal{D}_M$ is defined by \eqref{eq13}. In this domain the function $f(t,u,\bar{u},y,\bar{y},v,\bar{v},z,\bar{z})$ satisfies the Lipscchitz conditions in the variables $u,\bar{u},y,\bar{y},v,\bar{v},z,\bar{z}$ with the coefficients $L_0=1.30, L_1=7.20, L_2=0.47, L_3=0, L_4=0.94, L_5=0.32, L_6=0.31, L_7=0.86$. Therefore, $q=0.6446 <1$ and the problem \eqref{eqExam3} has a unique solution and the iterative method \eqref{iter1d}-\eqref{iter3d} converges. The results of the convergence of the problem are given in Table \ref{table3}.

\begin{table}[ht!]
\centering
\caption{The convergence in Example 3. }
\label{table3}
\begin{tabular}{ccccc}
\hline 
$N$ &	$h^2$ &	$K$	&$Error$&$Error1$  \\
\hline 
50	& 4.0000e-04 &	9&	9.2553e-08 & 3.0639e-07\\
100	& 1.0000e-04 &	9&	2.3182e-08 & 7.6420e-08\\
150	&4.4444e-05	&  9& 1.0307e-08 &  3.3946e-08 \\
200	&2.5000e-05	&11	&  5.7979e-09 & 1.9088e-08 \\
300	&1.1111e-05	& 9	&  2.5771e-09  & 8.4831e-09 \\
400&	6.2500e-06&	9	& 1.4497e-09& 4.7715e-09\\
500	&4.0000e-06	& 9	& 9.2781e-10& 3.0537e-09\\
800	&1.5625e-06&	8&	 3.6243e-10 &  1.1928e-09\\
1000&	1.0000e-06	&8	&2.3196e-10& 7.6339e-10\\
\hline 
\end{tabular} 
\end{table}
\noindent From Table \ref{table3}  it is easy to see that for this example when the right hand side is rather complicated, the errors of the approximate solution and its first derivative are not of order 4 as in the previous examples (Bica's examples), where  the right hand sides are simple. Obviously, here $Error$ and $Error1$ are of $O(h^3)$.

\begin{figure}[ht]
\begin{center}
\includegraphics[height=6cm,width=9cm]{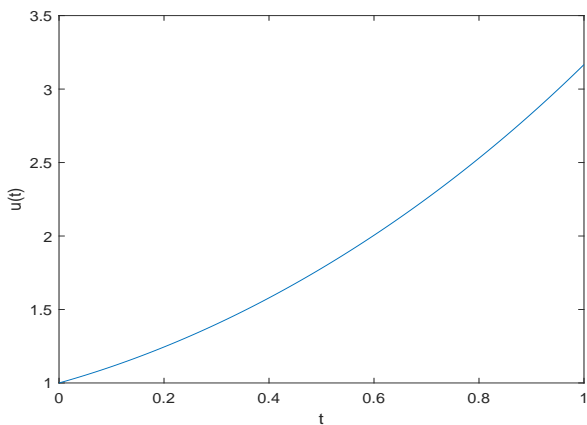}
\caption{The graph of the approximate solution in Example $4$. }
\label{fig1}
\end{center}
\end{figure}

\noindent {\bf Example 4 .} Consider the following problem
\begin{equation}\label{eqExam4}
\begin{split}
u^{(4)}(t)&= t^2-\frac{1}{4}u(t)+\frac{1}{4}[u(\varphi_0(t))]^2+\frac{1}{2} u'(t)u'(\varphi_1(t))+
           \frac{1}{8}[u''(t)+u''(\varphi_2(t))]u(t)\\
           &+\frac{1}{4}[\sin (u'''(t))+\cos (u'''(\varphi_3(t)))],\; t\in (0,1)\\
u(0)&=1, \ u(1)= \frac{19}{6}, \ u'(0)=1, \ u'(1)=\frac{7}{2},
\end{split}
\end{equation}
where $\varphi_0(t)=\frac{t}{2}, \varphi_1(t)= t^2, \varphi_2=\frac{t^2}{2}, \varphi_3(t)=\frac{t^2}{3}$.
At this moment we do not know any information of the solvability of the problem.\\
For the above problem
\begin{align*}
&f= t^2-\frac{1}{4}u +\frac{1}{4}\bar{u}^2+ y \bar{y}+\frac{1}{8}(v+\bar{v})u+\frac{1}{4}(\sin (z) 
+\cos (\bar{z})),\\
&g(t)=\frac{t^3}{6}+t^2+t+1.
\end{align*}
Hence, $\|g\|=\frac{19}{6}, \|g'\|=\frac{7}{2}, \|g''\|=3, \|g'''\|=1$. It is possible to verify that with $M=23$ there holds $|f| \le M$ in the domain $\mathcal{D}_M$,  where $\mathcal{D}_M$ is defined by \eqref{eq13}. Also, in 
$\mathcal{D}_M$ the function  $f(t,u,\bar{u},y,\bar{y},v,\bar{v},z,\bar{z})$ satisfies the Lipscchitz conditions in the variables $u,\bar{u},y,\bar{y},v,\bar{v},z,\bar{z}$ with the coefficients $L_0=1.48, L_1=1.62, L_2=L_3=3.7, L_4= L_5=0.41, L_6= L_7=0.25$. Therefore, $q=0.3857 <1$ and the problem \eqref{eqExam4} has a unique solution and the iterative method \eqref{iter1d}-\eqref{iter3d} converges. The results of the convergence of the problem is given in Table \ref{table4}. Recall that in the table $K$ is the number of iterations performed when 
 $\| U _k -U_{k-1}\|_{\bar{\omega_h}} \le 10^{-16}$.
\begin{table}[ht!]
\centering
\caption{The convergence in Example 4. }
\label{table4}
\begin{tabular}{cc}
\hline 
$N$  &	$K$	 \\
\hline 
50	 &	12 \\
100	 &	13\\
150		& 13  \\
200		&12	\\
300		&13	\\
400	&14	\\
500	&14	\\
800	&	13\\
1000	&13	
\\
\hline 
\end{tabular} 
\end{table}
The graph of the found approximate solution is depicted in Figure \ref{fig1}.\\

\noindent {\bf Remark.} 
It should be remarked that Theorem \ref{thm1} and Theorem \ref{thm3} give only sufficient conditions for the existence and uniqueness of solution, and the convergence of the iterative method for finding the solution. In the case if these conditions are not satisfied the iterative method also can be convergent. Below we give some examples for illustrating this remark.\\
\noindent {\bf Example 5 .} Consider the following problem
\begin{equation}\label{eqExam5}
\begin{split}
u^{(4)}(t)&= [u''( \dfrac{t}{4})]^4,\; 0<t<1,\\
u(0)&=1, \ u(1)= e, \ u'(0)=1, \ u'(1)=e,
\end{split}
\end{equation}
which has the exact solution  $u(t)=e^{t}$.\\
For this example $f=f(t,\bar{v})= (\bar{v})^4$. It is impossible to find $M>0$ so that $|f| \le M$ in the domain
$\mathcal{D}_M =\{(t, \bar{v}) \mid 0\le t  \le 1, |\bar{v}| \le 8-2e+\frac{M}{12}   \} $. Recall that $ 8-2e =\|g''\|$ (see Example 3). Therefore, the theorems \ref{thm1} and \ref{thm2} are not applicable. Regardless of this, the iterative method \eqref{iter1d}-\eqref{iter3d} converges to the exact solution  $u(t)=e^{t}$. Below are the results of convergence for the problem \eqref{eqExam5} (see Table \ref{table5}).
\begin{table}[ht!]
\centering
\caption{The convergence in Example 5. }
\label{table5}
\begin{tabular}{ccccc}
\hline 
$N$ &	$h^2$ &	$K$	&$Error$&$Error1$  \\
\hline 
50	& 4.0000e-04 &	20&	4.3691e-08 & 1.3415e-07\\
100	& 1.0000e-04 &	21&	1.0868e-08 & 3.3749e-08\\
150	&4.4444e-05	&  20& 4.8264e-09 &  1.5026e-08 \\
200	&2.5000e-05	&21	&  2.7140e-09 & 8.4562e-09 \\
300	&1.1111e-05	& 22	&  1.2059e-09  & 3.7597e-09 \\
400&	6.2500e-06&	22	& 6.7828e-10& 2.1150e-09\\
500	&4.0000e-06	& 22	& 4.3409e-10& 1.3537e-09\\
800	&1.5625e-06&	22&	 1.6956e-10 &  5.2882e-10\\
1000&	1.0000e-06	&22	&1.0852e-10& 3.3845e-10\\
\hline 
\end{tabular} 
\end{table}

\noindent {\bf Example 6 .} Consider the following problem
\begin{equation}\label{eqExam6}
\begin{split}
u^{(4)}(t)&= (u(t))^2 + (u''(\dfrac{t^2}{2}))^4,\; t\in (0,1)\\
u(0)&=1, \ u(1)= \frac{19}{6}, \ u'(0)=1, \ u'(1)=\frac{7}{2},
\end{split}
\end{equation}
We do not know any information of the solution of the problem. And it is easy to verify that the conditions of Theorem \ref{thm1} are not satisfied. Therefore, we are not ensured of the existence of solutions and the convergence of the iterative method. But the results of numerical experiments show that the iterative method \eqref{iter1d}- \eqref{iter3d} has a good convergence. The number of iterations needed for achieving the tolerance $\| U _k -U_{k-1}\|_{\bar{\omega_h}} \le 10^{-16}$  is 15 or 16.\\

\section{Conclusion}
In this paper we have considered the Dirichlet problem for fully fourth order nonlinear functional differential equation, where its right hand side contains all lower derivatives of proportional delay arguments. By reducing the problem to operator equation for the right hand side nonlinear function we have established the existence and uniqueness of solution under some easily verified conditions and constructed an efficient iterative method for finding the solution at both continuous and discrete levels. A total error of the actual numerical solution  which consists of the error of the iterative process and the error of discretization at each iteration is obtained.  Many examples with numerical experiments confirm the applicability and validity of the obtained theoretical results. To our best knowledge this is first time a boundary value problem for fully fourth order nonlinear functional differential equation is studied.\par  
It should be said that the proposed approach here to boundary value problems for nonlinear functional differential equations is  very effective because it simultaneously gives the existence results and the convergence of the solution method which is very easily realized.  
The approach can be applicable to a wide class of boundary value problems for nonlinear functional differential equations of arbitrary order and any linear boundary conditions.

\section*{Appendix} Derivatives in $t$ of the Green function $G(t,s)$:\\
\begin{equation*}
\begin{aligned}
G(t,s)=\dfrac{1}{6}\left\{\begin{array}{ll}
s^2(1-t)^2(3t-s-2ts), \quad 0\le s \le t \le 1,\\
t^2(1-s)^2,(3s-t-2ts) \quad 0\le t \le s \le 1,\\
\end{array}\right.
\end{aligned}
\end{equation*}
\begin{equation*}
\begin{aligned}
G_1(t,s)=\frac{\partial G(t,s)}{\partial t})=\left\{\begin{array}{ll}
- (s^2(2t - 2)(s - 3t + 2st))/6 - (s^2(2s - 3)(t - 1)^2)/6, \quad 0\le s \le t \le 1,\\
- (t^2(2s + 1)(s - 1)^2)/6 - (t(s - 1)^2(t - 3s + 2st))/3, \quad 0\le t \le s \le 1,\\
\end{array}\right.
\end{aligned}
\end{equation*}
\begin{equation*}
\begin{aligned}
G_2(t,s)=\frac{\partial ^2 G(t,s)}{\partial t^2})=\left\{\begin{array}{ll}
- (s^2(s - 3t + 2st))/3 - (s^2(2s - 3)(2t - 2))/3, \quad 0\le s \le t \le 1,\\
- ((s - 1)^2(t - 3s + 2st))/3 - (2t(2s + 1)(s - 1)^2)/3, \quad 0\le t \le s \le 1,\\
\end{array}\right.
\end{aligned}
\end{equation*}
\begin{equation*}
\begin{aligned}
G_3(t,s)=\frac{\partial ^3 G(t,s)}{\partial t^3})=\left\{\begin{array}{ll}
s^2(3-2s), \quad 0\le s < t \le 1,\\
s^2(3-2s)-1 \quad 0\le t < s \le 1,\\
\end{array}\right.
\end{aligned}
\end{equation*}
\begin{equation*}
\begin{aligned}
G_3^{*}(t,s)=\left\{\begin{array}{ll}
s^2(3-2s), \quad 0\le s < t \le 1,\\
t^2(3-2t)-1/2, \quad s=t \\ 
s^2(3-2s)-1 \quad 0\le t < s \le 1,\\
\end{array}\right.
\end{aligned}
\end{equation*}

\end{document}